\documentclass[reqno,final,12pt]{amsart}
\usepackage{amssymb,amsthm,amsmath,color,fullpage,url,thm-restate, booktabs}
\usepackage{comment}
\usepackage[noabbrev,capitalise]{cleveref}

\newtheorem{thm}{Theorem}
\newtheorem{lem}[thm]{Lemma}
\newtheorem{obs}[thm]{Observation}

\newtheorem{cor}[thm]{Corollary}

\newtheorem{fact}[thm]{Fact}
\newtheorem{rem}[thm]{Remark}
\newtheorem{claim}{Claim}[thm]

\newcommand{\eps}{\varepsilon}
 
\newcommand{\mc}[1]{\mathcal{#1}} 
\newcommand{\bb}[1]{\mathbb{#1}}

\newcommand{\bl}[1]{\boldsymbol{#1}}

\usepackage[noabbrev,capitalise]{cleveref}

\newcommand{\fs}[2]{\left(\frac{#1}{#2}\right)}
\newcommand{\s}[1]{\left(#1\right)}
\begin{document}	
\title{Optimizing the CGMS upper bound on Ramsey numbers}
\author{Parth Gupta}
\curraddr{School of Computer Science, University of Waterloo, Waterloo, ON, Canada.}		
\author{Ndiam\'e Ndiaye}
\curraddr{AGATE Center, Charles University, Prague, Czechia.}
\author{Sergey Norin}	
\author{Louis Wei}
\curraddr{University of Toronto, Toronto, ON, Canada}
\address{Department of Mathematics and Statistics, McGill University, Montr\'{e}al, QC, Canada.}
\thanks{PG, NN and SN were partially supported by an NSERC Discovery Grant. PG was also supported by an ISM Undergraduate Research Scholarship.
}
\begin{abstract} 
	In a recent breakthrough, Campos, Griffiths, Morris and Sahasrabudhe obtained the first exponential improvement of the upper bound on the diagonal Ramsey numbers since 1935. We shorten their proof, replacing the underlying \emph{book algorithm} with a simple inductive statement. This modification allows us to \begin{itemize} \item give a very short proof of an improved upper bound on the off-diagonal Ramsey numbers, which extends to the multicolor setting,
		\item clarify the dependence of the bounds on underlying parameters and optimize these parameters, obtaining, in particular, an upper bound
		$$R(k,k) \leq (3.8)^{k+o(k)}$$ 
		on the diagonal Ramsey numbers.
		\end{itemize}  
\end{abstract}
	\maketitle

\section{Introduction} The \emph{Ramsey number} $R(k, \ell)$ is the smallest positive integer $N$ such that
in any red-blue coloring of the edges of the complete graph on $N$ vertices there exists either a complete subgraph on $k$ vertices with all edges colored red (\emph{a red $K_k$}) or a complete subgraph on $\ell$ vertices with all edges colored blue (\emph{a blue $K_{\ell}$}).
Ramsey numbers were introduced by Ramsey ~\cite{R30} in 1930, and the problem of estimating their value has been one of the central questions in extremal combinatorics ever since, with particular attention paid to the \emph{diagonal} Ramsey numbers $R(k,k)$.
	
Erd\H{o}s and Szekeres ~\cite{ES35} proved that $R(k,k) \leq 4^k$ in 1935. Despite a series of important improvements ~\cite{C09,S20,T88} throughout the years, based on the quasirandomness properties of colorings close to the upper bound, the first exponential improvement of the Erd\H{o}s-Szekeres bound was obtained only very recently in a breakthrough result by Campos, Griffiths, Morris and Sahasrabudhe ~\cite{CGMS23}.
They proved that  there exists $\eps > 0$ such that 
\begin{equation}\label{e:CGMS}
R(k,k) \leq (4-\eps)^k
\end{equation} for sufficiently large $k \in \mathbb{N}$. While they mention that they have not attempted to seriously optimize  the value of $\eps$, focusing instead on a relatively simple proof, they show that \eqref{e:CGMS} is satisfied with $\eps = 2^{-7}$.   For the off-diagonal Ramsey numbers, Campos, Griffiths, Morris and Sahasrabudhe ~\cite{CGMS23} establish
\begin{equation}\label{e:CGMS2}
R(k, \ell) \leq e^{-\ell/400 + o(k)}\binom{k + \ell}{\ell}
\end{equation} for all $k, \ell \in \mathbb{N}$ with $\ell \leq k$, which again improves the best previously  known bound ~\cite{S20} by an exponential factor. 

The approach of~\cite{CGMS23} diverges from the previous methods and does not involve any notion of quasirandomness. Instead, it centres on the \emph{book algorithm}, which works with a pair of disjoint sets of vertices $(X,Y)$ with high density of red edges between them. Perhaps the most technical part of the algorithm is the way the changes in this density are controlled, regulated by a certain discrete logarithmic gradation of the interval $[p,1]$, where $p$ is approximately the initial density.

In an attempt to clarify the extent of applications of the approach of~\cite{CGMS23}, as well as the underlying optimization involved, we present a simpler proof of an improved upper bound on Ramsey numbers, based on a reinterpretation of their method. We do not introduce any new  combinatorial ideas. (In fact, we show that for $\ell < 0.69k$ one can abandon most of the ingredients used in~\cite{CGMS23}  and get an exponential improvement using a bare minimum.) We do, however, change the way we package and optimize applications of these ideas, replacing the CGMS book algorithm by a fairly simple induction, perhaps  closer paralleling the original Erd\H{o}s-Szekeres approach. In particular, we don't explicitly bound allowed changes in density of red edges, eliminating the technicality mentioned above.   

Using these modifications we obtain the following improved bounds. 

\vskip 5pt

In \cref{s:easy} we give a short proof of the following explicit bound on
Ramsey numbers
\begin{equation}\label{e:easy1}R(k,\ell) \leq  4(k+\ell)\fs{\left(\sqrt{5}+1\right) (k+2 \ell)}{4 \ell}^\ell \fs{k + 2\ell}{k}^{k/2}\end{equation}
for all positive integers $\ell \leq k$. (See \cref{c:easy}.)
The main idea of the proof is to inductively maintain a lower bound on the \emph{excess} number of red edges between the two sets of vertices $X$ and $Y$ under consideration, as compared to some fixed density $p$, i.e. a bound on  $e_R(X,Y)-p|X||Y|$ where $e_R(X,Y)$ denotes the number of red edges with one end in $X$ and another in $Y$.

It is not hard to verify that \eqref{e:easy1} gives an exponential improvement of the Erd\H{o}s-Szekeres bound for large enough $\ell$ and $k$, as long as $\ell \leq 0.69k$. It is meaningful even in the $\ell=o(k)$ regime, while~\cite{CGMS23} focuses on the regime $\ell=\Theta(k)$. Moreover, the results of \cref{s:easy} straightforwardly extend to the multicolor Ramsey numbers, as shown in \cref{s:multi}, while extending the whole book algorithm seems to require overcoming non-trivial technical obstacles.

\vskip 5pt
In \cref{s:book} we present our interpretation of the full power of the CGMS book algorithm, optimizing some of the parameters involved. We replace the quantity maintained inductively by a more involved one, which  essentially corresponds to a higher moment of the above mentioned excess number of edges. In \cref{s:opt} we optimize the initial value of density of red edges at which we can start our induction and prove our main result, the following. 

\begin{restatable}{thm}{Main}\label{t:main} 
	For all positive integers $\ell \leq k$	$$ R(k,\ell) \leq e^{G(\ell/k)k+o(k)}\binom{k+\ell}{\ell},  $$
	where $G(\lambda)= \s{-0.25\lambda + 0.03 \lambda^2 + 0.08\lambda^3}e^{-\lambda}.$
\end{restatable}

The term $o(k)$ in \cref{t:main} and in similar statements in this paper is uniform for  $1 \leq \ell \leq k$. That is \cref{t:main} states that there exists $\eta(k)$ with $\eta(k)/k \to 0$ as $k \to \infty$ such that $ R(k,\ell) \leq e^{G(\ell/k)k+
	\eta(k)}\binom{k+\ell}{\ell}  $ for all positive integers $\ell \leq k$.

In particular, \cref{t:main} implies that $$R(k,k) \leq e^{-0.14e^{-1}k+o(k)}\binom{2k}{k} = (4e^{-0.14e^{-1}})^{k+o(k)}=(3.7992\ldots)^{k+o(k)}, $$
and, more generally, that
$$R(k,\ell) \leq e^{-0.14e^{-1}\ell+o(k)}\binom{k+\ell}{\ell} \leq e^{-\ell/20+o(k)}\binom{k+\ell}{\ell},$$
significantly improving \eqref{e:CGMS} and \eqref{e:CGMS2}. 

\subsection*{AI Disclosure} The main results of this paper were obtained in Summer 2024 without any AI use. Numerical calculations used in derivation of \cref{t:main} from \cref{t:general} were incorrectly justified in the earliest public version. These were corrected in the PhD thesis~\cite{NdiamePhD} of the second named author.

The current, shorter derivation of \cref{t:main} from \cref{t:general} was obtained by 
OpenAI GPT-5.6 Sol based on the outline provided by the authors. The output was checked and edited by the authors, who take full responsibility for its correctness. OpenAI GPT-5.6 Sol  was also used for proofreading the paper. 

OpenAI Codex was used to formalize the main results of the paper (\cref{t:main} and \cref{c:easy}) in Lean 4. The formalization is available at \url{https://github.com/snorin239/RamseyLean}.

\subsection{Notation} As mentioned earlier we denote by $e_R(X,Y)$ the number of red edges with one end in a set  $X$  and another in $Y$. Given a color $C$ and a vertex $v$ of our graph, we denote by  $N_C(v)$  the set of vertices $u$ such that the edge $uv$ is colored in color $C$. (Throughout the bulk of the paper we work with two colors: red denoted by $R$, and blue denoted by $B$, but in \cref{s:multi} we introduce additional colors.) 

\section{Easy bound far from the diagonal}\label{s:easy}

In this section we present a very short inductive argument yielding an exponential improvement to the Erd\H{o}s-Szekeres upper bound on $R(k,\ell)$ for $\ell < 0.69k.$ 

The main object of our investigation is the same as in~\cite{CGMS23}: a pair of disjoint sets of vertices, the density of red edges between which is the critical parameter controlled during the inductive argument. We formalize the setting in the following definitions.
Let $X,Y$ be two non-empty disjoint subsets of vertices of a complete graph with edges colored red and blue. We say that $(X,Y)$ is a \emph{candidate}. We say that a candidate $(X,Y)$ is \emph{$(k,\ell,t)$-good} if $X \cup Y$ contains a red $K_k$ or $X$ contains a blue $K_t$ or $Y$ contains a blue $K_\ell$. The \emph{density of $(X,Y)$} is  $$d(X,Y) = \frac{e_R(X,Y)}{|X||Y|}.$$ For later convenience we extend this definition to pairs of not necessarily non-empty vertex sets by setting $d(X,Y)=0$ if $X = \emptyset$ or $Y=\emptyset$.
For $p \in [0,1]$ let $f_{p}(X,Y)=e_R(X,Y) - p|X||Y|$ denote the excess amount of red edges between $X$ and $Y$ when compared to density $p$. 
The following lemma shows that replacing $Y$  by $N_R(v) \cap Y$ for $v \in X$ on average  reduces $f_p$ by a factor at most $p$. It is essentially a restatement of~\cite[Observation 5.5]{CGMS23}.   

 \begin{lem}\label{l:FpAvg}
	Let $(X,Y)$ be a candidate and $p \in [0,1]$. Then \begin{equation}
	\sum_{v\in X}f_p(X,N_R(v) \cap Y) \geq p|X|f_p(X,Y).
	\end{equation}
	
\end{lem}
\begin{proof}
	We have
	\begin{align*}
	\sum_{v\in X}&f_p(X,N_R(v) \cap Y) - p|X|f_p(X,Y) \\&= \sum_{v \in X} \s{\sum_{y \in N_R(v) \cap Y}\s{ |N_R(y) \cap X| - p|X|} } -p|X|e_R(X,Y)+p^2|X|^2|Y| \\
	&= \sum_{y \in Y}\s{ 	\sum_{v \in N_R(y) \cap X} \s{|N_R(y) \cap X| -p|X|} -p|X||N_R(y) \cap X| +p^2|X|^2} \\&=\sum_{y \in Y} \s{ |N_R(y) \cap X| -p|X|}^2 \geq 0.
	\end{align*}
\end{proof}	

Next we need the following form of the Erd\H{o}s-Szekeres upper bound on Ramsey numbers.

\begin{obs}\label{o:easybound}For all $0<x<1$ and all positive integers $k$ and $\ell$
$$R(k,\ell) \leq x^{-k+1}(1-x)^{-\ell+1}.$$	
\end{obs}
\begin{proof}
	By induction on $k+\ell$. If $k=1$ or $\ell=1$ the statement clearly holds, and the induction step follows, as

	\begin{align*}R(k,\ell) &\leq R(k,\ell-1)+R(k-1,\ell) \leq x^{-k+1}(1-x)^{-(\ell-1)+1}+x^{-(k-1)+1}(1-x)^{-\ell+1} \\ &= x^{-k+1}(1-x)^{-\ell+1} (x+ (1-x)) = x^{-k+1}(1-x)^{-\ell+1} .
		\end{align*}
\end{proof}

The following lemma contains the main inductive argument we use to replace the book algorithm.

\begin{lem}\label{l:easy}
	Let $0<x<p<1$, let $k,\ell$ and $t$ be positive integers and let   $(X,Y)$ be a candidate such that 
	\begin{equation}\label{e:easy}
	f_{p}(X,Y) \geq  (k+t)x^{-k+1}(1-x)^{-\ell+1}(p-x)^{-t+1}
	\end{equation}	
	then $(X,Y)$ is $(k,\ell,t)$-good.
\end{lem}
\begin{proof}
	The proof is by induction on $k+t$. If $k=1$ or $t=1$ then  every candidate is $(k,\ell,t)$-good. This implies the base case of induction and allows us to assume $k,t \geq 2$ in the induction step. 
	
	By \cref{l:FpAvg} there exists $v \in X$ such that  $f_p(X,N_R(v) \cap Y) \geq p \cdot f_p(X,Y)$. Let $Y' = N_R(v) \cap Y, X_B=  N_B(v) \cap X$ and $X_R=  N_R(v) \cap X$. We will prove that $(X,Y)$ is $(k,\ell,t)$-good by showing that either $(X_R,Y')$ is $(k-1,\ell,t)$-good or 
	$(X_B,Y')$ is $(k,\ell,t-1)$-good by the induction hypothesis, or $|Y|$ is large enough to contain a red $K_k$ or a blue $K_{\ell}$ by \cref{o:easybound}.
	
	If $$f_p(X_R,Y') \geq \frac{k+t-1}{k+t} x f_p(X,Y) \geq (k+t-1)x^{-(k-1)+1}(1-x)^{-\ell+1}(p-x)^{-t+1 }$$ then by the induction hypothesis $(X_R,Y')$ is  $(k-1,\ell,t)$-good, i.e. $X_R \cup Y'$ contains a red $K_{k-1}$ or $X_R$ contains a blue $K_t$ or $Y'$ contains a blue $K_\ell$. In the first case,  $X_R \cup Y' \cup \{v\} \subseteq X \cup Y$ contains a red $K_k$, as $X_R \cup Y' \subseteq N_R(v).$ It follows that in each case  $(X,Y)$ is  $(k,\ell,t)$-good, as desired, and so we may assume that $f_p(X_R,Y') < \frac{k+t-1}{k+t} x f_p(X,Y)$.
	
Symmetrically, we may assume that $f_p(X_B,Y') <\frac{k+t-1}{k+t}(p-x) f_p(X,Y)$, as otherwise $(X_B,Y')$ is  $(k,\ell,t-1)$-good by the induction hypothesis, implying that  $(X,Y)$ is  $(k,\ell,t)$-good. (If $X_B$ contains a blue $K_{t-1}$ then $X_B \cup \{v\}$ contains a blue $K_t$).

It follows that
\begin{align*}
 	 p  f_p(X,Y) &\leq f_p(X,Y') = f_p(X_R,Y') + f_p(X_B,Y') + f_p(\{v\},Y') \\ &< \frac{k+t-1}{k+t}pf_p(X,Y) + |Y|,
	\end{align*}
and so $|Y| \geq \frac{1}{{k+t}}pf_p(X,Y)$.
Thus, using \eqref{e:easy} to lower bound $f_p(X,Y)$, we have
$$
|Y| \geq  px^{-k+1}(1-x)^{-\ell+1}(p-x)^{-t+1} \geq x^{-k+1}(1-x)^{-\ell+1} \geq  R(k,\ell),
$$	
where the last inequality uses \cref{o:easybound}, implying that $(X,Y)$ is  $(k,\ell,t)$-good. 
\end{proof}	

Finally, we derive from~\cref{l:easy} the promised bound on the Ramsey numbers by induction on $\ell$, applying~\cref{l:easy}, if the density of red edges is sufficiently high, and, otherwise, applying the induction hypothesis to the blue neighborhood of an arbitrary vertex.  

\begin{thm}\label{t:easy}
	For all $\frac{\sqrt{5}-1}{\sqrt{5}+1}< p <1$ and all positive integers $k$ and  $\ell$ $$R(k,\ell) \leq 4(k+\ell)\s{\frac{1+\sqrt{5}}{2}p+\frac{1-\sqrt{5}}{2}}^{-k/2}(1-p)^{-\ell} .$$
\end{thm}
\begin{proof}
	By induction on $\ell$. Note that the theorem trivially holds if $k=1$ or $\ell =1 $. In particular, the base case $\ell=1$ trivially holds and we assume $k, \ell \geq 2$ in the induction step.
	
	Let $x = \frac{1+\sqrt{5}}{2}p+\frac{1-\sqrt{5}}{2} >0$. Then $0<x<p<1$ and note that 	\begin{equation}\label{e:PEquation}
(1-x)(p-x) = \s{\frac{1+\sqrt{5}}{2}-\frac{1+\sqrt{5}}{2}p}\s{\frac{\sqrt{5}-1}{2}-\frac{\sqrt{5}-1}{2}p}=(1-p)^2.
	\end{equation} Consider a red-blue coloring of edges of a complete graph on $n \geq \lfloor 4(k+\ell)x^{-k/2}(1-p)^{-\ell} \rfloor$ vertices. 
	In particular, we have \begin{align*} p(1-p)^2n \geq  p(1-p)^2 \s{4(k+ \ell)x^{-1} (1-p)^{-2} - 1}   \geq  4(k+\ell) - p \geq 4(k+\ell)(1-p),\end{align*}
	where the second inequality uses $x \leq p$.
	Thus 
	\begin{equation}\label{e:nLower} 
		n \geq \frac{4(k + \ell)}{p(1-p)}.
	\end{equation} 
	 If $|N_B(v)|\geq  \lfloor 4(k+\ell-1)x^{-k/2}(1-p)^{-\ell+1} \rfloor$ for some vertex $v$ then the coloring of $N_B(v)$ contains a red $K_k$ or a blue $K_{\ell-1}$ by the induction hypothesis, and so our coloring contains a red $K_k$ or a blue $K_\ell$. 
	 
	Thus we may assume that $$|N_B(v)| \leq   \lfloor 4(k+\ell-1)x^{-k/2}(1-p)^{-\ell+1} \rfloor - 1  \leq 4(k+\ell-1)x^{-k/2}(1-p)^{-\ell+1} - 1$$ for every $v$, i.e. \begin{align*} |N_R(v)|  &= n-1 - |N_B(v)| \geq n - 1 -  \s{4(k+\ell-1)x^{-k/2}(1-p)^{-\ell+1} -1} \\&> n - \frac{k+\ell-1}{k+\ell}(1-p)(n+1) > \s{p +\frac{(1-p)}{k+\ell}}n - 1. \end{align*} Choose a random partition $(X,Y)$ by assigning every vertex inependently to $X$ or $Y$, each with probability $1/2$. Then
	$$\bb{E}[e_R(X,Y)] = \frac{1}{4}\sum_{v}|N_R(v)| >  \s{p +\frac{(1-p)}{k+\ell}}\frac{n^2}{4} -\frac{n}{4}.$$ 
	It follows that 
	\begin{align*}\bb{E}[f_p(X,Y)] &\geq \bb{E}[e_R(X,Y)] -\frac{pn^2}{4} \geq \frac{(1-p)}{k+\ell}\frac{n^2}{4} -\frac{n}{4}\\& \geq \frac{(1-p)^2}{k+\ell}\frac{n^2}{4}	\geq (k+\ell)x^{-k}(1-p)^{-2\ell+2}\\&=  (k+\ell)x^{-k}(1-x)^{-\ell+1}(p-x)^{-\ell+1},
	\end{align*}
	where the second inequality holds by \eqref{e:nLower}, and the last equality uses \eqref{e:PEquation}. Thus there exists a candidate $(X,Y)$ for which $f_p(X,Y) \geq (k+\ell)x^{-k}(1-x)^{-\ell+1}(p-x)^{-\ell+1}$. By \cref{l:easy} $(X,Y)$ is $(k,\ell,\ell)$-good, implying that  our  coloring contains a red $K_k$ or a blue $K_\ell$, as desired.
\end{proof}	

Substituting $$p=\frac{(\sqrt{5}+1) k+(2\sqrt{5}-2)\ell}{\left(\sqrt{5}+1\right) (k+2 \ell)}$$ in \cref{t:easy} yields the following corollary, mentioned in the introduction.

\begin{cor}\label{c:easy}
	For all positive integers $k \geq \ell$ $$R(k,\ell) \leq 4(k+\ell)\fs{\left(\sqrt{5}+1\right) (k+2 \ell)}{4 \ell}^\ell \fs{k + 2\ell}{k}^{k/2}. $$
\end{cor}

Let $ES(k,\ell)=\binom{k+\ell-2}{k-1}$ denote the Erd\H{o}s-Szekeres upper bound on the Ramsey numbers. Then $ES(k,\ell)=e^{O(\log k)}\fs{k+\ell}{k}^{k}\fs{k+\ell}{\ell}^{\ell},$  and so \cref{c:easy} implies
\begin{align*} \frac{R(k,\ell)}{ES(k,\ell)} &\leq e^{O(\log k)}\fs{(\sqrt{5}+1)(k+2\ell)}{4(k+\ell)}^{\ell}\fs{(k+2\ell)k}{(k+\ell)^2}^{k/2}. 
\end{align*}
Thus \cref{c:easy} yields an improvement  exponential in  $\ell$ of the Erd\H{o}s-Szekeres bound for  $\log k \ll \ell < 0.6989k,$ and for $\ell=o(k)$ the improvement is of the order $e^{O(\log k)}\fs{\sqrt{5}+1}{4}^{(1+o(1))\ell} < e^{-0.21\ell + 	 O(\log k)}$.
For comparison, the strongest bound for moderately small $\ell$ given in \cite{CGMS23} is of the form $$R(k,\ell) \leq e^{-0.05\frac{k}{k+\ell}\cdot\ell+o(k)}ES(k,\ell)$$
for $\ell \leq k/9$ (see ~\cite[Theorem 9.1]{CGMS23}).
Meanwhile, in the same regime ($\ell \leq k/9$) \cref{c:easy} implies 
$$R(k,\ell) \leq e^{-0.16\ell+ O(\log k)}ES(k,\ell).$$
In the following section, using the full power of the book algorithm, we further improve on this bound and extend the improvement to the diagonal case.

\section{Optimizing the book algorithm}\label{s:book}	

In this section we present extensions of \cref{l:easy} and \cref{t:easy},  which we use to give an upper bound $R(k,\ell)$ that improves on the bounds from ~\cite{CGMS23}. Just like in \cref{l:easy} the main graph-theoretic ideas are borrowed from~\cite{CGMS23}, but now we need to refine our inductive statement and to use every trick in the book from~\cite{CGMS23}, making the argument substantially more technical. Namely:
\begin{itemize}
	\item instead of $e_R(X,Y) - p|X||Y|$ we lower bound ``higher moments'' of density $$(e_R(X,Y) - p|X||Y|)^r|X|^{1-r}|Y|^{1-r}$$ in the regime $r \to \infty$,
	\item we implement the arguments corresponding to \emph{big blue steps} and \emph{degree regularization} steps of the~\cite{CGMS23} book algorithm, which we did not need in~\cref{l:easy}.  
\end{itemize}	

Instead of \cref{l:FpAvg}   we will need the following variant of ``the convexity of density'' bound.
 
 \begin{lem}\label{l:FpAvg2}
 	Let $(X,Y)$ be a candidate. Then \begin{equation}
 	\sum_{v\in X}d(X,N_R(v) \cap Y)|N_R(v) \cap Y| \geq e_R(X,Y)d(X,Y).
 	\end{equation}
 \end{lem}
 \begin{proof}
We have
 	\begin{align*}
 |X|&\s{\sum_{v\in X}d(X,N_R(v) \cap Y)|N_R(v) \cap Y|} = \sum_{v \in X}e_R(X,N_R(v) \cap Y) \\ & = \sum_{v\in X} \s{\sum_{y \in N_R(v) \cap Y} |N_R(y) \cap X|} = \sum_{y \in Y}|N_R(y) \cap X|^2 \\ &\geq  |Y|\fs{\sum_{y \in Y}|N_R(y) \cap X|}{|Y|}^2 =\frac{(e_R(X,Y))^2}{|Y|}  =|X| e_R(X,Y)d(X,Y),
 	\end{align*}
 	where the inequality holds by convexity of the square function.
 \end{proof}	

We say that a pair $(S,T)$ of vertex-disjoint sets in a complete graph with edges colored red and blue is \emph{a blue book} if all the edges with both ends in $S$ and all edges with one end in $S$ and another in $T$ are colored blue. Blue books are useful for finding blue cliques
as a blue $K_{a}$ in $T$ can be extended to a blue $K_{a+|S|}$  by adding the vertices of $S$ to it. 

The next lemma allows us to extract a large blue book from $X$ if we find sufficiently many vertices with large blue neighborhood. It is essentially ~\cite[Lemma 4.1]{CGMS23} with a different choice of parameters. The proof is exactly the same, but we reproduce it for completeness. To do so we need the following inequality relating binomial coefficients used in the proof of~\cite[Lemma 4.1]{CGMS23}. (In ~\cite{CGMS23} the inequality below is stated only for $0<\sigma < 1$, but it trivially holds for $\sigma =1$.)

\begin{fact}{~\cite[Fact 4.2]{CGMS23}}\label{f:binomial} Let $0< \sigma \leq 1$, let $b,m$ be positive integers with $b \leq \sigma m /2$. Then
\begin{equation}\label{e:binomial}
	\binom{\sigma m}{b} \geq \sigma^b\binom{m}{b}\exp\s{-\frac{b^2}{\sigma m}}
\end{equation}
\end{fact}	

\begin{lem}\label{l:BBook}
	Let $0< \mu < 1$, let $b,k,m$ be positive integers  with $m \geq 10\mu^{-1}b^2$. Let $(X,Y)$ be a candidate such that $|X| \geq 5m^2$, and there exist at least $R(k,m)$ vertices $v \in X$ such that 
	$$ |N_B(v) \cap X| \geq \mu|X|.$$
	Then $X$ contains a red $K_k$ or a blue book $(S,T)$ with $|S| \geq b$ and $|T| \geq  \frac{\mu^b}{2}|X|$.
\end{lem}

\begin{proof}
	Let $W$ be the set of vertices $v \in X$ such that 	$ |N_B(v) \cap X| \geq \mu|X|.$ As $|W| \geq R(k,m)$, $W$ contains a  red $K_k$ or a blue $K_m$. In the first case the lemma holds, so we assume the second, and  let $U$ be the set of vertices of the blue $K_m$  in $W$. Let $\sigma = \frac{e_B(U,X\setminus U)}{|U||X \setminus U|}$. Then  \begin{equation}\label{e:sigma} \sigma = \frac{ \sum_{v \in U} |N_B(v) \cap (X \setminus U)|}{|U||X \setminus U|} \geq  \frac{\mu|X|-m}{|X|}  \geq \mu \s{1- \frac{1}{5b}} \geq \frac{1}{2}\mu,\end{equation}
	where the second to last inequality holds as $|X| \geq 5m^2 \geq  5\mu^{-1}mb$.
	Also note that \begin{equation}\label{e:sigmaAvg}\sum_{v \in X \setminus U}|N_B(v) \cap U| =  e_B(U,X\setminus U) = \sigma m |X \setminus U|,\end{equation}
	by definition of $\sigma$.
	
	Let $S$ be a subset of $U$ of size $b$ chosen uniformly at random, and let $T$ be the set of all vertices  $v \in X \setminus U$ such that $S \subseteq N_B(v)$. Then $(S,T)$ is a blue book, and it suffices to show that $|T| \geq  \frac{\mu^b}{2}|X|$ for some choice of $S$.  We have \begin{align*}\bb{E}[|T|] &=\binom{m}{b}^{-1}\sum_{v \in X \setminus U}\binom{|N_B(v) \cap U|}{b} \\ &\stackrel{\eqref{e:sigmaAvg}}{\geq} \binom{m}{b}^{-1}\binom{\sigma m}{b}|X \setminus U|  &\mathrm{by \;convexity\; of}\; x \to \binom{x}{b} \\ & \stackrel{\eqref{e:binomial}}{\geq} \sigma^b\exp\s{-\frac{b^2}{\sigma m}}|X \setminus U| \\ &\geq \sigma^b e^{-1/5}|X \setminus U| & \mathrm{as}\qquad \frac{b^2}{\sigma m} \leq \frac{b^2}{\sigma \cdot 10 \mu^{-1}b^2} \leq \frac{1}{5} \\ &\stackrel{\eqref{e:sigma}}{\geq}  \mu^b\s{1- \frac{1}{5b}}^b  e^{-1/5}\s{1-\frac{m}{|X|}}|X| \\ &\geq \mu^b \cdot \frac{4}{5} \cdot e^{-1/5} \cdot \frac{4}{5}|X|  & \mathrm{as}\qquad \s{1- \frac{1}{5b}}^b  \geq \frac{4}{5}  \\ &\geq \frac{\mu^b}{2}|X|, \end{align*}
implying the existence of the desired choice of $S$.
\end{proof}	

In \cref{s:easy} we used \cref{o:easybound} to upper bound Ramsey numbers in what was essentially the base case of the book algorithm, or, more precisely, the case when we consider a candidate $(X,Y)$ with $|X|$ small. In this section we tighten all aspects of our argument and, in  particular, instead of \cref{o:easybound} we would like to apply increasingly stronger upper bounds of the same form, which we iteratively obtain.

We use the following slightly technical setup to encode the bounds we can use. 
Let $\mc{R}_0$ be the set of all pairs $(x,y) \in (0,1)^2$ such that there exists $N=N(x,y)$ such that \begin{equation}\label{e:R} R(k,\ell) \leq x^{-k}y^{-\ell}  \qquad \mathrm{for \; all} \; k, \ell \in \bb{N} \; \mathrm{such\; that \;} k+\ell \geq N .\end{equation} Let $\mc{R}$ be the closure of $\mc{R}_0$, and let $\mc{R_*}$ be the interior of $\mc{R}$.

The next easy observation records the properties of $\mc{R}$ and $\mc{R_*}$.

\begin{obs}\label{o:r}
	\begin{enumerate}
		\item $(x,1-x) \in \mc{R}$ for all $0\leq x \leq 1$,
		\item if  $(x,y) \in \mc{R}$, $0 < x' \leq x, 0 < y' \leq y$  then $(x',y') \in \mc{R}$,
		\item if  $(x,y) \in \mc{R}$, $0 < x' < x, 0 < y'< y$  then $(x',y') \in \mc{R}_*$,
		\item if $R(k,\ell) \leq x^{-k-o(k)}y^{-\ell-o(\ell)}$ for some $(x,y) \in (0,1)^2$ then $(x,y) \in \mc{R}$.
	\end{enumerate}	
\end{obs}	

\begin{proof}
\cref{o:easybound} implies (1). As $ x^{-k}y^{-\ell} \leq (x')^{-k}(y')^{-\ell}$ for $0 < x' \leq x, 0 < y' \leq y$,  (2) holds with $\mc{R}$ replaced by $\mc{R}_0$, which in turn implies  that it holds for  $\mc{R}$.  The item (2) in turn implies (3). 

Finally, we verify (4). Let $(x,y) \in (0,1)^2$ be such that $R(k,\ell) \leq x^{-k-o(k)}y^{-\ell-o(\ell)}$. It suffices to show that for every $ \eps > 0$ we have $(x^{1+\eps},y^{1+\eps}) \in \mc{R}_0$, i. e. that there exists $N$ such that \begin{equation}\label{e:4}R(k,\ell) \leq x^{-k(1+\eps)}y^{-\ell(1+\eps)}\end{equation} for all positive integers $k$ and $\ell$ such that $k + \ell \geq N$. Let $n$ be such that $$R(k,\ell) \leq x^{-k(1+\eps/2)} y^{-\ell(1+\eps/2)}$$ for all positive integers $k,\ell \geq n$.
Let $N \geq 2n$ be chosen so that $$\max\{x^{(N-n)\eps/2},y^{(N-n)\eps/2}\} \leq (xy)^{n(1+\eps/2)}$$ and let  $k$ and $\ell$ such that $k + \ell \geq N$. If $k, \ell \geq n$ then \eqref{e:4} holds as desired, therefore we assume without loss of generality that $\ell \leq n$ and so $k \geq N-n$. Then \begin{align*} 
R(k,\ell) &\leq R(k,n) \leq x^{-k(1+\eps/2)}y^{-n(1+\eps/2)} = x^{-k(1+\eps)} \cdot x^{k\eps/2}y^{-n(1+\eps/2)} \\& \leq x^{-k(1+\eps)}   \cdot x^{(N-n)\eps/2}(xy)^{-n(1+\eps/2)} \leq x^{-k(1+\eps)} \leq  x^{-k(1+\eps)}y^{-\ell(1+\eps)},
\end{align*}
as desired.
\end{proof} 

As noted in the above proof the Erd\H{o}s-Szekeres argument via \cref{o:easybound} implies \cref{o:r} (1). It is not hard to see that conversely \cref{o:r} implies the Erd\H{o}s-Szekeres upper bound on Ramsey numbers up to subexponential factors. The problem of extending the set of points which provably belong to $\mc{R}$ can be considered as  dual to the problem of obtaining asymptotic upper bounds on Ramsey numbers $R(k, \ell)$ in the regime $k = \Theta(\ell)$. While it is not clear why the optimal upper bounds on Ramsey numbers must come in the form of \eqref{e:R}, the bounds in this form are very convenient to apply inductively, as we've already done in \cref{l:easy} and will do again in \cref{t:bookmain}.
Eventually, we will show that \begin{equation}\label{e:R2} (x,\frac{1}{4}e^{ 0.14e^{-1}}x^{-1}) \in \mc{R} \qquad \mathrm{for\: } \qquad 0.4663 \leq x  \leq 0.564\end{equation} (see \cref{r:final}). Substituting $x=0.51304\ldots$, the solution of $x = \frac{1}{4}e^{0.14e^{-1}}x^{-1}$  we obtain $R(k,k) \leq x^{-2k+o(k)} = (3.7992\ldots)^{k+o(k)}$, which is exactly the bound obtained from~\cref{t:main}.

Note that \cref{t:easy} implies that $$\s{\s{\frac{1+\sqrt{5}}{2}p+\frac{1-\sqrt{5}}{2}}^{1/2},1-p} \in \mc{R} \qquad \mathrm{for} \qquad \frac{\sqrt{5}-1}{\sqrt{5}+1}< p <1.$$

Finally, we need the following straightforward limit evaluation.

\begin{lem}\label{l:limit} 
	For all $0 <  p ,  \mu < 1$ 
	\begin{equation}\label{e:limit} \lim_{r \to \infty}(p^{1/r}-\mu)^r(1-\mu)^{1-r}=p^{\frac{1}{1-\mu}}(1-\mu)
	\end{equation} 
\end{lem}

\begin{proof}
	We have
	\begin{align*} \lim_{r \to \infty}&\log\s{(p^{1/r}-\mu)^r(1-\mu)^{-r}} \\&=\lim_{r \to \infty} \s{r \log\s{1+\frac{p^{1/r}-1}{1-\mu}}} =\frac{\lim_{r \to \infty} \s{r(p^{1/r}-1)}}{1-\mu}=\frac{\log p}{1-\mu},
	\end{align*}
	implying \eqref{e:limit}.
\end{proof}	

We are now ready for the main technical result of this section, which tightens~\cref{l:easy}.

\begin{lem}\label{t:bookmain}  For all $0< \mu_0,x_0,y_0,p <1$ such that $x_0 < p^{\frac{1}{1-\mu_0}}(1-\mu_0)$ and $(x_0,y_0) \in \mc{R_*}$ there exists $L_0$ such that for all positive integers $k,\ell,t$ with $\ell \geq L_0$  the following holds.
			Let $(X,Y)$ be a candidate such that $d(X,Y) \geq p$ and  \begin{equation}\label{e:bookmain} |X||Y| \geq x_0^{-k}y_0^{-\ell}\mu_0^{-t }\end{equation} 
			then $(X,Y)$ is $(k,\ell,t)$-good.
\end{lem}

\begin{proof} First, let us briefly outline the proof. We follow the same general strategy as the proof of \cref{l:easy}: we  find a vertex $v \in X$ such that $d(X, N_R(v) \cap Y)$ is not much smaller than $d(X,Y)$ with the goal of applying the induction hypothesis to the candidate $(N_R(v) \cap X,N_R(v) \cap Y)$ or  $(N_B(v) \cap X, N_R(v) \cap Y)$. There are however several important adjustments.
\begin{itemize}
	\item The induction hypothesis \eqref{e:moment} involves a higher moment of excess as mentioned in the beginning of the section.
	\item We start by proving a lower bound of roughly $p|Y|$ on $|N_R(v) \cap Y|$ for every $v \in X$. (If this bound does not hold we apply the induction hypothesis to $(X \setminus \{v\},Y)$). A similar argument could have been used in \cref{l:easy}, but was not needed.
	\item More importantly, we show that we can assume that  $|N_R(v) \cap X|$ is not much larger than $\mu|X|$ for most of  the vertices of $X$, as otherwise, using \cref{l:BBook} we can find a large blue book $(S,T)$ in $X$ and apply the induction hypothesis to $(T,Y)$.
	\item Finally, using these degree bounds and convexity of density (\cref{l:FpAvg2}), we find suitable $v \in X$. 
\end{itemize}

We start by quantifying the ``extra room'' implicit in the strict inequality $x_0 < p^{\frac{1}{1-\mu_0}}(1-\mu_0)$ and the condition $(x_0,y_0) \in \mc{R_*}$. More precisely, by \cref{l:limit}  and the condition $x_0 < p^{\frac{1}{1-\mu_0}}(1-\mu_0)$, there exists $r>1$ such that $p^{1/r} > \mu_0$ and $x_0  <  (p^{1/r} -  \mu_0)^{r}(1-\mu_0)^{1-r}$. This condition further implies that  $x_0+\mu_0<1$. Thus we deduce that   there exists  $\eps > 0$ such that $$ (1+\eps)(\mu_0+\eps) \leq 1, \quad \mu_0+x_0+2\eps \leq 1, \quad (x_0+3\eps, y_0 +3\eps) \in \mc{R}, \quad (p-\eps)^{1/r}  \geq  \mu_0 +3\eps,  $$
 \begin{equation}\label{e:xx} \mu_0 + 2\eps \leq \frac{\mu_0}{\mu_0 + 2\eps +x_0},   \quad  p \geq 2\eps, \quad \mathrm{and} \quad x_0  \leq  ((p-\eps)^{1/r} -  \mu_0-3\eps)^{r}(1-\mu_0-\eps)^{1-r} -\eps.\end{equation} 

We choose $L_0$ implicitly to be sufficiently large as a function of $\eps$ and $r$ to satisfy the  inequalities throughout the proof.
 Let $x=x_0 +\eps$, $y=y_0+\eps, \mu = \mu_0+\eps$,	$\delta_{n}=\frac{\eps}{n}$. Note that by \eqref{e:xx} we have $x <1$ as  $$ x = x_0 + \eps < 1- \mu_0 - \eps \leq \frac{1}{1+\eps}.$$ Thus $x_0 = x-\eps \leq \frac{x}{1+\eps}$. Similarly, $y_0  \leq \frac{y}{1+\eps}$ and  $\mu_0 \leq \frac{\mu}{1+\eps}$. 
 
We will prove by induction on $k+ t$ for fixed $\ell \geq L_0$  that if $(X,Y)$ is a candidate with $d(X,Y) \geq p -\delta_{k+t}$ such that \begin{equation}\label{e:moment}
	(d(X,Y) + \delta_{k+t} - p)^r |X||Y| \geq x^{-k}y^{-\ell}\mu^{-t }
\end{equation}  then 
$(X,Y)$ is $(k,\ell,t)$-good. 

Note that
\begin{align*} (d(X,Y) &+ \delta_{k+t} - p)^r |X||Y| \stackrel{\eqref{e:bookmain}}{\geq} 
	 \frac{\eps^r}{(k+t)^r}x_0^{-k}y_0^{-\ell}\mu_0^{-t } \notag\\ &\geq \frac{\eps^r}{(k+t)^r}(1+\eps)^{k+\ell+t} x^{-k}y^{-\ell}\mu^{-t } \geq x^{-k}y^{-\ell}\mu^{-t },
\end{align*}
where the last inequality holds as long as $L_0$ is sufficiently large as a function of $\eps$ and $r$.
Thus the above statement implies the lemma. 
If $k=1$ or $t=1$ then $(X,Y)$ is trivially $(k,\ell,t)$-good, implying, in particular the base case of our statement. Thus we move on to the induction step and assume $k,t \geq 2$.
 
We assume without loss of generality that $|N_R(v) \cap Y| \geq (p-\delta_{k+t})|Y|$ for every $v \in X$, as otherwise we can replace $X$ by $X \setminus \{v\}$ increasing the left side of \eqref{e:moment} as  $$(d(X,Y) + \delta_{k+t} - p)^r |X||Y| = (d(X,Y) + \delta_{k+t} - p)^{r-1}(e_R(X,Y)-(p-\delta_{k+t})|X||Y|),$$
and both terms on the right side of this identity increase after the replacement. Thus   \begin{equation}\label{e:degreg}d(X',Y) \geq p-\delta_{k+t}\end{equation}  for every $X' \subseteq X, X' \neq \emptyset.$ For readers familiar with~\cite{CGMS23} we will indicate  how the steps in our proof  correspond to steps of  their book algorithm. The observation in this paragraph corresponds to the \textbf{degree regularization step} of the book algorithm.

Next we use \eqref{e:moment} to obtain a lower bound on $|X|$ by first upper bounding $|Y|$. 
If $|Y| \geq (x+\eps)^{-k} (y + \eps)^{-\ell}$ then $Y$ contains a red $K_k$ or a blue $K_{\ell}$  as $(x+\eps, y +\eps) \in \mc{R_*}$ by \cref{o:r} (3) and so  $(X,Y)$ is $(k,\ell,t)$-good.
Thus we may assume that $|Y| \leq  (x+\eps)^{-k} (y + \eps)^{-\ell}$ and so \begin{align} \label{e:x}|X| &\stackrel{\eqref{e:moment}}{\geq}  \frac{ x^{-k}y^{-\ell}\mu^{-t }}{(d(X,Y) + \delta_{k+t} - p)^r|Y|} \geq \frac{ x^{-k}y^{-\ell}\mu^{-t }}{ (x+\eps)^{-k} (y +\eps)^{-\ell}} \notag\\& \geq \fs{x+\eps}{x}^k\fs{y+\eps}{y}^{\ell}\mu^{-t} \geq (1+\eps)^{k +\ell+t}.\end{align}

Our next goal is to show that most of the vertices of $X$ have blue degree at most $(\mu+\eps)|X|$ as otherwise using \cref{l:BBook}  we can apply the induction hypothesis to the candidate $(T,Y)$ where $(S,T)$ is a large blue book in $X$ guaranteed by 
\cref{l:BBook} and \eqref{e:x}. This part of the proof corresponds to the \textbf{big blue step}.

Let $$b =\left\lceil \frac{2r \log(k+t)- r\log{\eps} + \log{2}}{\log(1+\eps)}\right \rceil,  m = \lceil  10\mu^{-1}b^2 \rceil, \qquad \mathrm{and} \qquad w=R(k,m) $$
with the first two being the parameters we will use in  \cref{l:BBook}.
Note that $m \leq C\log^2(k+t),$ where $C$ is a constant depending only on $r$ and $\eps$.  As $R(k,m) \leq k^m \leq \exp(C\log^3(k+t))$, by \eqref{e:x} we can choose $L_0$ large enough so that \begin{equation}\label{e:x2}|X| \geq 5m^2 \qquad \mathrm{and} \qquad
w \leq \frac{\eps(p-\eps)}{(k+t)^3}|X|.\end{equation}

Let  $W =\{ v \in X \: | \:  |N_B(v) \cap X| \geq (\mu + \eps)|X|\}$. Suppose first that $|W| \geq w$. Then by \cref{l:BBook}, $X$ contains a red $K_k$ or a blue book $(S,T)$ with $|S| \geq b$ and $|T| \geq  \frac{ (\mu + \eps)^b}{2}|X|$. In the first case $(X,Y)$ is $(k,\ell,t)$-good, and thus we may assume that the second case holds.  If $b \geq t$ then $X$ contains a blue $K_t$ and so  $(X,Y)$ is once again $(k,\ell,t)$-good, and so we further assume $b < t$.
As $d(T,Y) \geq p-\delta_{k+t}$ by \eqref{e:degreg}, we have $d(T,Y) + \delta_{k+t-b} -p \geq \delta_{k+t-1}-\delta_{k+t} \geq \frac{\eps}{(k+t)^2}.$
Thus \begin{align*} 
(d(T,Y) + \delta_{k+t-b} - p)^r |T||Y| &\geq \fs{\eps }{(k+t)^2}^r \frac{(\mu + \eps)^b}{2}|X||Y| \\ &\stackrel{\eqref{e:moment}}{\geq}  \fs{\eps }{(k+t)^2}^r\frac{(\mu + \eps)^b}{2}x^{-k}y^{-\ell}\mu^{-t } \\& = \frac{1}{2}\fs{\eps}{(k+t)^2}^r\fs{\mu + \eps}{\mu}^b x^{-k}y^{-\ell}\mu^{-t+b} \\&\geq x^{-k}y^{-\ell}\mu^{-t+b},
\end{align*}
where the last inequality holds by the choice of $b$ as
$$ \fs{\mu + \eps}{\mu}^b \geq (1+\eps)^b \geq 2\fs{(k+t)^{2}}{\eps}^r.$$ 
Thus $(T,Y)$ is $(k,\ell,t-b)$-good by the induction hypothesis, implying that   $(X,Y)$ is $(k,\ell,t)$-good. 

Thus we may assume that $|W| \leq w$. We now move on to the last part of the argument, corresponding to either the \textbf{red step} or the \textbf{density increment step} of the book algorithm, where the steps are symmetric to each other in our implementation. 

Let $p'= p-\delta_{k+t-1} \geq p-\eps.$ By \cref{l:FpAvg2}, we have \begin{equation}
\sum_{v\in X}d(X,N_R(v) \cap Y)|N_R(v) \cap Y| \geq e_R(X,Y)d(X,Y).
\end{equation}

As $$\sum_{v\in W}d(X,N_R(v) \cap Y)|N_R(v) \cap Y| \leq w|Y| \stackrel{\eqref{e:x2}}{\leq}\frac{\eps}{(k+t)^3}(p -\eps)|X||Y|  \leq  \frac{\eps}{(k+t)^3} e_R(X,Y)$$
we have  \begin{equation*}
\sum_{v\in X - W}d(X,N_R(v) \cap Y)|N_R(v) \cap Y| \geq \s{d(X,Y)-\frac{\eps}{(k+t)^3}}e_R(X,Y).
\end{equation*}

Thus there exists $v\in X$ such that $|N_B(v) \cap X| \leq (\mu+\eps)|X|$ and  \begin{equation}\label{e:alpha}d(X,N_R(v) \cap Y) \geq d(X,Y)-\frac{\eps}{(k+t)^3}.\end{equation}  

Let \begin{align*} &X_R = N_{R}(v) \cap X,\qquad &X_B =N_{B}(v) \cap X,\qquad  &Y'=N_{R}(v) \cap Y, \\ &\alpha = d(X,Y')-p', \qquad &\alpha_R = d(X_R,Y')-p', \qquad &\alpha_B = d(X_B,Y')-p'.\end{align*} Note that $|Y'| \geq (p-\eps)|Y|$ by \eqref{e:degreg},
 \begin{align*} \alpha &= d(X,Y')-p' \stackrel{\eqref{e:alpha}}{\geq} d(X,Y) -\frac{\eps}{(k+t)^3}+ \delta_{k+t-1} - p \\ &\geq  \delta_{k+t-1} - \delta_{k+t}  -\frac{\eps}{(k+t)^3} =  \frac{\eps}{(k+t)(k+t-1)}-\frac{\eps}{(k+t)^3} \geq \frac{\eps}{(k+t)^2},\end{align*}  and 
\begin{align*} 
(\alpha_R&|X_R|+\alpha_B|X_B| + 1)|Y'|\\ &=(e_R(X_R,Y')+e_R(X_B,Y')+e_R(\{v\},Y'))-p'(|X_R|+|X_B|)|Y'|\\ &\geq e_R(X,Y')-p'|X||Y'| = \alpha|X||Y'|,
\end{align*}
implying
\begin{equation}\label{e:moment2}
\frac{\alpha_R}{\alpha}\frac{|X_R|}{|X|} +\frac{\alpha_B}{\alpha}\frac{|X_B|}{|X|} + \frac{1}{\alpha|X|}\geq 1.
\end{equation} 
If $\alpha_R \geq 0$ and 
$\alpha_R^r|X_R| \geq \frac{x\alpha^r|X|}{p-\eps}$ then \begin{align*} (d(X_R,Y')-p')^r|X_R||Y'| &=\alpha_R^r|X_R||Y'|  \geq x(d(X,Y')  - p')^r|X|\frac{|Y'|}{p-\eps} \\&\stackrel{\eqref{e:alpha}}{\geq} x\s{d(X,Y) -\frac{\eps}{(k+t)^3}+ \delta_{k+t-1} - p}^r|X||Y| \\ &\geq x(d(X,Y) + \delta_{k+t} - p)^r|X||Y| {\geq} x^{-k+1}y^{-\ell}\mu^{-t}.\end{align*}
This implies that $(X_R,Y')$ is $(k-1,\ell,t)$-good by the induction hypothesis, and so $(X,Y)$ is $(k,\ell,t)$-good.

Thus we may assume that either $\alpha_R < 0$ or  $\alpha_R^r|X_R| < \frac{x\alpha^r|X|}{p-\eps}$, i.e. $$\frac{\alpha_R}{\alpha} \frac{|X_R|}{|X|}< x^{1/r}(p-\eps)^{-1/r}\fs{|X_R|}{|X|}^{1-1/r}.$$ Symmetrically, we also have $$\frac{\alpha_B}{\alpha} \frac{|X_B|}{|X|}<  \mu^{1/r}(p-\eps)^{-1/r}\fs{|X_B|}{|X|}^{1-1/r}.$$
It follows from \eqref{e:moment2} that
\begin{equation}\label{e:moment3}
x^{1/r}\fs{|X_R|}{|X|}^{1-1/r}+\mu^{1/r}\fs{|X_B|}{|X|}^{1-1/r}+ \frac{(p-\eps)^{1/r}}{\alpha|X|} >(p-\eps)^{1/r}. 
\end{equation}

The remainder of the proof is occupied with showing that \eqref{e:moment3} contradicts \eqref{e:xx}.
As $\frac{|X_R|}{|X|}+\frac{|X_B|}{|X|}\leq 1$, $\frac{|X_B|}{|X|} \leq (\mu+\eps)$, the function $\lambda \to x^{1/r}(1-\lambda)^{1-1/r} + \mu^{1/r}\lambda^{1-1/r} $ increases for $0 \leq \lambda \leq \frac{\mu}{\mu+x}$ and $ \mu+\eps \leq \frac{\mu}{\mu+x}$ we can upper bound the left side of \eqref{e:moment3} by replacing $\frac{|X_B|}{|X|} $ by $\mu+\eps$  and $\frac{|X_R|}{|X|} $ by $1-\mu$, implying
\begin{equation}\label{e:moment4}
x^{1/r}(1-\mu)^{1-1/r}+\mu+ \eps + \frac{(p-\eps)^{1/r}}{\alpha|X|} > (p-\eps)^{1/r}. 
\end{equation}
As $\alpha \geq \frac{\eps}{(k+t)^2}$ we have
$$ \frac{(p-\eps)^{1/r}}{\alpha|X|}  \stackrel{\eqref{e:x}}{\leq} \frac{(k+t)^2}{\eps(1+\eps)^{k+\ell+t}} \leq \eps,$$
where the last inequality holds whenever $L_0$ (and thus $\ell$) is large enough as a function of $\eps$.
It follows that \eqref{e:moment4} in turn implies
$x^{1/r}(1-\mu)^{1-1/r}+\mu+ 2\eps  > (p-\eps)^{1/r}, $ i.e.
$$ x > ((p-\eps)^{1/r} - \mu- 2\eps)^r(1-\mu)^{1-r} $$
contradicting \eqref{e:xx}.
\end{proof}

The main result of this section follows immediately from \cref{t:bookmain}.

\begin{thm}\label{t:bookCor}
For all $0< \mu,x,y,p <1$ such that $x < p^{\frac{1}{1-\mu}}(1-\mu)$ and $(x,y) \in \mc{R_*}$  there exists $L_0$ such that for all positive integers $k,\ell$ with $\ell \geq L_0$ the following holds. Every red-blue coloring of the edges of the complete graph on $N \geq x^{-k/2}(\mu y)^{-\ell/2}$ with the density of red edges  at least $p$ contains a red $K_k$ or a blue $K_{\ell}$. 
\end{thm}
\begin{proof}
Let $y_0 > y$ such that  $(x,y_0) \in \mc{R_*}$  and let $L_0$ be chosen sufficiently large to satisfy the requirements below including the condition $(y_0/y)^{\ell/2} \geq 3$ for $\ell\geq L_0$. Then $N  \geq x^{-k/2}(\mu y)^{-\ell/2} \geq 2x^{-k/2}(\mu y_0)^{-\ell/2} + 1$, and so $\lfloor N/2 \rfloor \geq x^{-k/2}(\mu y_0)^{-\ell/2}$. Let $X$ and $Y$ be two disjoint subsets of vertices of our graph each of size $\lfloor N/2 \rfloor$ chosen to maximize $d(X,Y)$. Then $d(X,Y) \geq p$ and $|X||Y| = (\lfloor N/2 \rfloor)^2 \geq x^{-k}y_0^{-\ell}\mu^{-\ell}$. The theorem follows by applying \cref{t:bookmain} to the candidate $(X,Y)$ with $x_0=x,\mu_0=\mu$.   
\end{proof}

\section{Optimizing descent to a candidate}\label{s:opt}

In this section we derive our main result from \cref{t:bookCor} in the same manner that \cref{t:easy} was derived from \cref{l:easy} in \cref{s:easy}: By using induction on $\ell$ with the induction hypothesis applied to the blue neighborhood of a vertex, if the density of red edges is insufficiently high to apply  \cref{t:bookCor}.

We express our upper bound on Ramsey numbers $R(k,\ell)$ with $\ell \leq k$ in the form $e^{F(\ell/k)k~+~o(k)}$ and seek to minimize $F$. As $e^{F((\ell-1)/k)k}/e^{F(\ell/k)k} \approx e^{-F'(\ell/k)},$ we can use the induction hypothesis  if the density of blue edges is about $e^{-F'(\ell/k)}$. Thus we will be applying  \cref{t:bookCor} when the density of red edges is at least $1 - e^{-F'(\ell/k)}$. The parameters $x, y$ and $\mu$ are then chosen as functions of $\ell/k$ to satisfy the requirements in \cref{t:bookCor}. We formalize this strategy in the following theorem.

\begin{thm}\label{t:general}
	Let $F:(0,1] \to \bb{R}_+$ be smooth and let $M,X,Y:(0,1] \to (0,1)$ with $M$ continuous be such that  
	\begin{align*} &F'(\lambda) > 0, \qquad  X(\lambda) = (1 - e^{-F'(\lambda)})^{\frac{1}{1-M(\lambda)}}(1-M(\lambda)), \qquad (X(\lambda),Y(\lambda)) \in \mc{R},
\\ & and  \qquad F(\lambda) > -\frac{1}{2} \s{\log X(\lambda)+ \lambda \log M(\lambda) + \lambda\log Y(\lambda)} \end{align*}
for all $0< \lambda \leq 1$. 
Then $$R(k,\ell) \leq e^{F(\ell/k)k + o(k)},$$	
	for all $k \geq \ell$.
\end{thm}

\begin{proof} We need to show that for every $\eps > 0$ we have \begin{equation}\label{e:epsBound}R(k,\ell) \leq e^{(F(\ell/k)+\eps)k}\end{equation} for all $k$ sufficiently large as a function of $\eps,F,M,X$ and $Y$ and $\ell \leq k$. As $R(k,\ell) = e^{o(k)}$ for  $\ell=o(k)$, it suffices to establish \eqref{e:epsBound} for $\ell > \eps'k$ for some $\eps'>0$ depending only on $\eps$.

Our goal is to prove the theorem by induction on $\ell$. We need, however, to resolve a technical issue arising from the fact that \cref{t:bookCor} is only applicable for $\ell$ large enough as a function of parameters of the theorem, which we have not made explicit. We obtain a uniform lower bound on the size of $\ell$ for any given $\eps$ by applying  \cref{t:bookCor} for finitely many choices of parameters  in the following claim.

\begin{claim}\label{c:gen}
	There exist $\delta, L >0$  such that for all positive integers  $k \geq \ell \geq \max\{ L, \eps' k \}$ the following holds.  Every red-blue coloring of the edges the complete graph on $N \geq e^{(F(\ell/k)+\eps)k}$ with the density of red edges at least $(1 - e^{-F'(\ell/k)+\delta})$ contains a red $K_k$ or a blue $K_{\ell}$. 
\end{claim}

\begin{proof}  
	Let $\delta > 0$ be such that $$e^{-\delta}\fs{1 - e^{-F'(\lambda)+2\delta}}{1-e^{-F'(\lambda)}}^{\frac{1}{1-M(\lambda)}} \geq e^{-\eps + \delta} \quad \mathrm{and} \quad F'(\lambda) > 2 \delta$$ for all $\lambda \in [\eps',1]$, and let $\eps''$ be such that $|F'(\lambda)-F'(\lambda')| \leq \delta$ and $|F(\lambda)-F(\lambda')| \leq \eps/2$ for all $\lambda,\lambda' \in [\eps',1]$ such that $|\lambda-\lambda'| \leq \eps''$.
	
	Let $\Lambda=\{1 - i\eps'' \: |\: i=0,1,\ldots, \lfloor (1-\eps')/\eps''\rfloor\} $.
	By  \cref{t:bookCor} there exists $L$ such that for every $\lambda \in \Lambda$ the outcome of \cref{t:bookCor}  applied with parameters $p_\lambda=1 - e^{-F'(\lambda)+2\delta}$,  $\mu_{\lambda}=M(\lambda)$, $x_{\lambda} = e^{-\delta}(p_{\lambda})^{\frac{1}{1-\mu_{\lambda}}}(1-\mu_{\lambda})$ and $y_{\lambda} = e^{-\delta}Y(\lambda)$ holds for all  positive integers $k,\ell$ with $\ell \geq L.$ Note that \begin{align*} X(\lambda) \geq x_{\lambda} = e^{-\delta}X(\lambda)\fs{1 - e^{-F'(\lambda)+2\delta}}{1-e^{-F'(\lambda)}}^{\frac{1}{1-\mu_{\lambda}}} \geq e^{-\eps+\delta}X(\lambda). \end{align*} 
	
	Given  $k \geq \ell \geq \max\{ L, \eps' k \}$ let $\lambda \in \Lambda$ be such that $\ell/k \leq \lambda \leq \ell/k + \eps''$. The density of red edges in the coloring we consider is at least $ 1 - e^{-F'(\ell/k)+\delta} \geq 1 - e^{-F'(\lambda)+2\delta}$ by the choice of $\delta$. Thus by the choice of $L$  the claim holds by \cref{t:bookCor}  applied with parameters corresponding to $\lambda$ as \begin{align*} N &\geq \exp\s{(F(\ell/k)+\eps)k} \geq \exp\s{(F(\lambda)+\eps/2)k} \\ &\geq    \exp\s{-\frac{1}{2} \s{k\log X(\lambda)+ \lambda k \log M(\lambda) + \lambda k \log Y(\lambda)}+\eps k/2} \\& \geq (e^{-\eps}X(\lambda))^{-k/2}(\mu_{\lambda} e^{\delta} y_{\lambda})^{-\ell/2} \geq (e^{-\eps +\delta}X(\lambda))^{-k/2}(\mu_{\lambda} y_{\lambda})^{-\ell/2} \\ &\geq x_{\lambda}^{-k/2}(\mu_{\lambda}  y_{\lambda})^{-\ell/2}.\end{align*} 	
\end{proof}

With \cref{c:gen} in hand it is not hard to finish the proof by induction on $\ell$. The base case $\ell \leq \eps'k$ holds as noted above.
For the induction step consider a red-blue coloring of the edges the complete graph on $N \geq e^{(F(\ell/k)+\eps)k}$ vertices. If the density of red edges is at least $(1 - e^{-F'(\ell/k)+\delta})$ then the coloring contains a red $K_k$ or a blue $K_{\ell}$ by \cref{c:gen}. Otherwise there exists a vertex $v$ with $$\deg_B(v) \geq e^{-F'(\ell/k)+\delta}N -1 \geq e^{-F'(\ell/k)+\delta/2}N  \geq \exp\s{F(\ell/k)k+\eps k-F'(\ell/k)+\delta/2},$$ 
where the second inequality holds  as $N \geq e^{\eps  k} \geq \frac{e^{F'(\ell/k)}}{e^{\delta/2}-1}$ for sufficiently large $k$.
Thus it suffices to show that $\exp\s{F(\ell/k)k+\eps k-F'(\ell/k)+\delta/2} \geq R(k,\ell-1)$. By the induction hypothesis, this is implied by \begin{equation}\label{e:gen2}F(\ell/k)k+\delta/2-F'(\ell/k) \geq F\s{\frac{\ell-1}{k}}k.\end{equation} As $F$ is smooth there exists $\lambda \in \left[\frac{\ell-1}{k},\frac{\ell}{k}\right]$ such that
$k(F(\ell/k)-F((\ell-1)/k)) = F'(\lambda)$, and so \eqref{e:gen2} can be rewritten as $\delta/2 \geq F'(\ell/k) -F'(\lambda).$ This last inequality holds for sufficiently large $k$, as $|F'(\ell/k) -F'(\lambda)| \leq \frac{1}{k}\max_{x \in [\eps'/2,1]}|F''(x)| \leq \delta/2.$
\end{proof}

We now approximate the optimal parameters to use in \cref{t:general} in stages. These were obtained by first finding continuous, piece-wise linear $F(\lambda)$ and piece-wise constant $M(\lambda)$ so that $F$ and $M$ satisfy the conditions of \cref{t:general} with $X(\lambda)= (1 - e^{-F'(\lambda)})^{\frac{1}{1-M(\lambda)}}(1-M(\lambda))$ and $Y(\lambda)=1-X(\lambda)$. (Except, of course, the condition that $F$ is smooth.)
Such $F$ and $M$  are built iteratively, respectively linear and constant, on intervals $[i/N,(i+1)/N]$ for $i=0,1,\ldots, N-1$, where we used $N=10$. The resulting $F$ and $M$ are approximated by reasonably simple smooth functions, and $M(\lambda)= \lambda  e^{-\lambda}$ turns out to be a reasonably good choice and, in fact, the behavior of $F$ is not too sensitive to small changes in $M$. This is already sufficient to obtain a bound $R(k,k) \leq (3.87)^{k+o(k)}$.

To improve this bound further we use the following lemma, which leverages our improved upper bounds on $R(k, \ell)$, allowing us to increase $Y(\lambda)$, while maintaining the condition $(X(\lambda),Y(\lambda)) \in \mc{R}$.  
For the remainder of this section let
\begin{align*}
	h(\lambda)&=(1+\lambda)\log(1+\lambda)-\lambda\log\lambda,\\
	g_b(\lambda)&=\bigl(-\tfrac14\lambda+b\lambda^2
	+0.08\lambda^3\bigr)e^{-\lambda},\\
	F_b(\lambda)&=h(\lambda)+g_b(\lambda)
	\qquad (0<\lambda\leq 1).
\end{align*}

\begin{lem}\label{lem:frontier}
	Let $f\colon(0,1]\to\mathbb R$ be differentiable and strictly concave,
	with $f'(t)>0$, and suppose that
	\begin{equation}\label{eq:f-bound}
		R(k,\ell)\leq \exp\bigl(f(\ell/k)k+o(k)\bigr)
		\qquad (1\leq \ell\leq k).
	\end{equation}
	Define
	\[
	A(t)=e^{-f'(t)},\qquad B(t)=e^{tf'(t)-f(t)}.
	\]
	Suppose that $A(t)<B(t)$ for $0<t\leq1$, that $A$ increases from $0$ to
	$a=A(1)$, and that $B$ decreases from $1$ to $b=B(1)$.  Then
	\begin{equation}\label{eq:frontier}
		Y_f(x)=
		\begin{cases}
			B(t),&0<x<a,\quad A(t)=x,\\[1mm]
			e^{-f(1)}/x,&a\leq x\leq b,\\[1mm]
			A(t),&b<x<1,\quad B(t)=x,
		\end{cases}
	\end{equation}
	is well defined, takes values in $(0,1)$, and
	\[
	(x,Y_f(x))\in\mc{R}
	\qquad (0<x<1).
	\]
\end{lem}

\begin{proof}Fix \(t\in(0,1]\), and abbreviate \(A=A(t)\) and \(B=B(t)\).
	By concavity, for every \(0<s\leq1\),
	\begin{equation*}
		f(s)+\log B+s\log A
		=f(s)-f(t)-(s-t)f'(t)\leq0.
	\end{equation*}
	Since \(A<B\), we have \(\log A<\log B\), and hence
	\[
	f(s)+\log A+s\log B
	=\bigl(f(s)+\log B+s\log A\bigr)
	+(1-s)(\log A-\log B)
	\leq0.
	\]
	Now let \(1\leq\ell\leq k\) and put \(s=\ell/k\).  Combining these
	two inequalities with \eqref{eq:f-bound} gives
	\[
	R(k,\ell)
	\leq B^{-k}A^{-\ell}e^{o(k)}
	\qquad\text{and}\qquad
	R(k,\ell)
	\leq A^{-k}B^{-\ell}e^{o(k)}.
	\]
	If \(1\leq k<\ell\), applying the same estimates to
	\(R(\ell,k)=R(k,\ell)\) instead gives
	\[
	R(k,\ell)
	\leq A^{-k}B^{-\ell}e^{o(\ell)}
	\qquad\text{and}\qquad
	R(k,\ell)
	\leq B^{-k}A^{-\ell}e^{o(\ell)}.
	\]  Applying
	\cref{o:r}(4) to the two estimates yields
	\[
	(A(t),B(t))\in\mc{R}
	\qquad\text{and}\qquad
	(B(t),A(t))\in\mc{R}.
	\]
	It remains to justify the middle part of \eqref{eq:frontier}.  From
	\[
	\log a=-f'(1)
	\qquad\text{and}\qquad
	\log b=f'(1)-f(1)
	\]
	we obtain \(ab=e^{-f(1)}\).  Fix \(a\leq x\leq b\), and set
	\[
	y=\frac{e^{-f(1)}}{x}=\frac{ab}{x}.
	\]
	In particular, \(a\leq y\leq b\), so \(y\in(0,1)\).
	
	By concavity  we have $
	f(s)\leq f(1)-(1-s)f'(1)	$ for $s \in (0,1]$.

	Using \(\log y=-f(1)-\log x\), we consequently have
	\begin{align*}
		f(s)+\log x+s\log y
		&=f(s)-sf(1)+(1-s)\log x\\
		&\leq(1-s)\bigl(f(1)-f'(1)+\log x\bigr)\\
		&=(1-s)(\log x-\log b)\leq0,
	\end{align*}
	where the last inequality uses \(x\leq b\).  Similarly,
	\begin{align*}
		f(s)+s\log x+\log y
		&=f(s)-f(1)-(1-s)\log x\\
		&\leq-(1-s)\bigl(f'(1)+\log x\bigr)\\
		&=(1-s)(\log a-\log x)\leq0,
	\end{align*}
	where the last inequality uses \(x\geq a\).
	
	If \(1\leq\ell\leq k\), applying the first inequality with
	\(s=\ell/k\) and using \eqref{eq:f-bound} gives
	\[
	R(k,\ell)
	\leq \exp\bigl(f(\ell/k)k+o(k)\bigr)
	\leq x^{-k}y^{-\ell}e^{o(k)}.
	\] 
	If \(1\leq k<\ell\), we obtain the same result using the second  inequality instead.
	Thus, 
	\[
	R(k,\ell)\leq x^{-k-o(k)}y^{-\ell-o(\ell)},
	\]and hence
	\((x,y)\in\mc{R}\) by \cref{o:r}(4).
	
	Finally, the hypotheses on \(A\) and \(B\) give existence and
	uniqueness of the parameter \(t\) in each of the two outer parts of
	\eqref{eq:frontier}.  At the junction points,
	\[
	\frac{e^{-f(1)}}{a}=\frac{ab}{a}=b=B(1),
	\qquad
	\frac{e^{-f(1)}}{b}=\frac{ab}{b}=a=A(1),
	\]
	so the three parts agree at their endpoints.  Consequently,
	\(Y_f\) is well defined and takes values in \((0,1)\).
\end{proof}

We record the elementary numerical verification used below.  It is
stated with substantially less precision than the computation provides.

\begin{lem}\label{lem:numerics}
	For $b\in\{0.08,0.03\}$ and a function $M\colon(0,1]\to(0,1)$, set
	\begin{equation}\label{eq:X-def}
		P_b(\lambda)=1-e^{-F_b'(\lambda)},\qquad
		X_{b,M}(\lambda)=(1-M(\lambda))
		P_b(\lambda)^{1/(1-M(\lambda))}.
	\end{equation}
	The following assertions hold for every $0<\lambda\leq1$.
	
	\begin{enumerate}
		\item Let
		\[
		M_0(\lambda)=\lambda e^{-\lambda},\quad
		X_0=X_{0.08,M_0},\quad Y_0=1-X_0.
		\]
		Then
		\begin{equation}\label{eq:base-slack}
			F_{0.08}(\lambda)+\frac12\bigl(
			\log X_0(\lambda)+\lambda\log M_0(\lambda)
			+\lambda\log Y_0(\lambda)\bigr)
			>2.4\cdot10^{-3}\lambda.
		\end{equation}
		
		\item Let
		\[
		M_1(\lambda)=\lambda
		\exp\!\left(-\frac23\lambda-\frac3{10}\lambda^2\right),
		\quad X_1=X_{0.03,M_1},\quad
		Y_1(\lambda)=Y_{F_{0.08}}(X_1(\lambda)),
		\]
		where $Y_{F_{0.08}}$ is defined in \cref{lem:frontier}.
		Then
		\begin{equation}\label{eq:final-slack}
			F_{0.03}(\lambda)+\frac12\bigl(
			\log X_1(\lambda)+\lambda\log M_1(\lambda)
			+\lambda\log Y_1(\lambda)\bigr)
			>3.0\cdot10^{-4}\lambda.
		\end{equation}
	\end{enumerate}
\end{lem}

\begin{proof}
	These are one-variable inequalities involving only elementary
	functions. The details of the numerical verification of these inequalities are given in Appendix A. The claimed numerical bounds are certified by Mathematica.
\end{proof} 
 
\Main*

\begin{proof} First note that, for
	$0.03\leq b\leq0.08$,
	\begin{equation}\label{eq:concavity}
		F_b''(t)=-\frac1{t(1+t)}+e^{-t}\bigl(
		0.5+2b+(0.23-4b)t+(b-0.48)t^2+0.08t^3\bigr)<0.
	\end{equation}
	Indeed, the polynomial in parentheses is at most $0.66$ for $b \in  [0.03, 0.08], t \in (0,1]$, while
	$0.66t(1+t)e^{-t} \leq 1.32/e<1$.  Also
	\begin{equation}\label{eq:positive}
		F_b(t)>t(\log2-\tfrac14)>0,
		\qquad F_b'(t)\geq\log2-\tfrac14>0.
	\end{equation}
	For the second inequality, if
	$q_b(t)=-t/4+bt^2+0.08t^3$, then
	$q_b'(t)-q_b(t)\geq-1/4$ for $b \in  [0.03, 0.08], t \in (0,1]$.
	
	Apply \cref{t:general} with $F=F_{0.08}$ and the functions $M_0,X_0,Y_0$
	from Lemma~\ref{lem:numerics}.  The identity required there for $X$
	is exactly \eqref{eq:X-def}; the point $(X_0 (\lambda),Y_0 (\lambda))$ belongs to $\mc{R}$ by
	the Erd\H{o}s--Szekeres bound \cref{o:easybound} and \cref{o:r} (4); and the strict numerical inequality is
	\eqref{eq:base-slack}.  Equations \eqref{eq:positive} and
	\eqref{eq:X-def} also show that all four functions take values in the
	ranges required by \cref{t:general}.  We obtain
	\begin{equation}\label{eq:base-bound}
		R(k,\ell)\leq
		\exp\bigl(F_{0.08}(\ell/k)k+o(k)\bigr)
		\qquad(1\leq\ell\leq k).
	\end{equation}
	
	We next verify that Lemma~\ref{lem:frontier} applies to
	$f=F_{0.08}$.  Strict concavity follows from \eqref{eq:concavity}, and
	\[
	A'(t)=-F_{0.08}''(t)A(t)>0,\qquad
	B'(t)=tF_{0.08}''(t)B(t)<0.
	\]
	Moreover, $A(t)\to0$ and $B(t)\to1$ as $t\downarrow0$.  Finally,
	\[
	\log\frac{B(t)}{A(t)}=(1+t)F_{0.08}'(t)-F_{0.08}(t)
	\]
	is decreasing, and its value at $t=1$ is $0.57/e>0$.  Thus all the
	hypotheses of Lemma~\ref{lem:frontier} hold.  By
	\eqref{eq:base-bound},
	\[
	(X_1(\lambda),Y_1(\lambda))\in\mc{R}
	\qquad(0<\lambda\leq1).
	\]
	
	Apply \cref{t:general} once more, now with
	$F=F_{0.03}$ and $M_1,X_1,Y_1$.  We have $0<M_1(\lambda)<1$, and
	\eqref{eq:positive} and \eqref{eq:X-def} imply
	$0<X_1(\lambda)<1$.  Lemma~\ref{lem:frontier} gives
	$0<Y_1(\lambda)<1$ and membership in $\mc{R}$.  The remaining strict
	inequality required by \cref{t:general} is precisely
	\eqref{eq:final-slack}.  Consequently,
	\[
	R(k,\ell)\leq
	\exp\bigl(F_{0.03}(\ell/k)k+o(k)\bigr).
	\]
	Uniformly for $1\leq\ell\leq k$, Stirling's formula gives
	\[
	\log\binom{k+\ell}{\ell}
	=h(\ell/k)k+o(k).
	\]
	Since $F_{0.03}=h+g_{0.03}$, this is exactly
	\[
	R(k,\ell)\leq
	e^{G(\ell/k)k+o(k)}\binom{k+\ell}{\ell},
	\qquad
	G(\lambda)=\bigl(-0.25\lambda+0.03\lambda^2
	+0.08\lambda^3\bigr)e^{-\lambda},
	\]
	as desired.
\end{proof}

\begin{rem}\label{r:final}
	Note that the calculations done for $F_{0.08}$ in the above proof of \cref{t:main} can be repeated for $F_{0.03}$ and \cref{lem:frontier} implies  \begin{equation}\label{e:R3}
(x, Y_{F_{0.03}}(x)) \in \mc{R} \qquad \mathrm{for\: every} \quad x \in (0,1). 
	\end{equation} This implies an earlier claim \eqref{e:R2}. 
	
	We asked ChatGPT 5.6 Sol to perform an additional iteration of optimization in Theorem 14. A preliminary, unverified iteration suggests that  \cref{t:main} holds with 
	$$G_{AI}(\lambda) = e^{-\lambda}(-0.3864\lambda+0.8347\lambda^2
	-2.0156\lambda^3+2.7171\lambda^4-1.7541\lambda^5
	+0.4522\lambda^6).$$
	If verified it would improve the upper bound on the diagonal Ramsey numbers to  $R(k,k) \leq (3.78233\ldots)^{k+o(k)}.$
	 
	 Further improvements by performing additional iterations are possible, but we expect that lowering the base of the exponent below $3.7$  and, likely, even below $3.75$ would require new ideas.
\end{rem}


\section{Multicolor Ramsey numbers}\label{s:multi}

In this section we  present extensions of the results in \cref{s:easy}  to multicolor Ramsey numbers. The techniques translate straightforwardly to the multicolor setting. In fact, the main reason that we present the bounds for two colors separately is to avoid complicating the notation.   

It will be convenient for us to use the following additional notation to encode the parameters of multicolor Ramsey numbers. Fix a positive integer $c$ for the remainder of this section. Let $\bl{x} \in \bb{R}^c$ be a vector. (We will use bold letters for vectors.) Then we denote its components by $x_1,\ldots,x_c$ and use the convention $x =\sum_{i=1}^cx_i$. Similarly $\bl{\ell}=(\ell_1,\ldots,\ell_c),$ and $\ell = \sum_{i=1}^c\ell_i$, etc. 

We will always distinguish a single color red ($R$) in the colorings we consider and  label the rest as $(B_1,\dots,B_c)$. For $k \in \bb{N}$ and  $\bl{\ell}\in \bb{N}^c$
\emph{the multicolor Ramsey number} $R(k,\bl{\ell})$ is the minimum integer $N$ such that every coloring of the edges of the complete graph on $N$ vertices in $c+1$ colors $R,B_1,\dots,B_c$ contains a complete subgraph on $k$ vertices with all edges colored $R$ or a complete subgraph on $\ell_i$ vertices with all edges colored $B_i$ for some $1 \leq i \leq c$.

The Erd\H{o}s-Szekeres argument to bound classical two-color Ramsey numbers naturally extends to multicolor ones, giving an upper bound on $R(k,\bl{\ell})$ of the form $$ES(k,\bl{\ell}) = e^{o(k+\ell)}\frac{(k+\ell)^{k+\ell}}{k^k \ell_1^{\ell_1}\cdots\ell_c^{\ell_c}} =  e^{o(k+\ell)} ES(k,\ell)\cdot \frac{\ell^{\ell}}{\ell_1^{\ell_1}\cdots\ell_c^{\ell_c}}.$$
The main change between our results in the two-color and multicolor cases is the introduction of the same factor  $\frac{\ell^{\ell}}{\ell_1^{\ell_1}\cdots\ell_c^{\ell_c}}$, which we denote by $\Theta(\bl{\ell})$.

\cref{o:easybound} extends as follows.

\begin{obs}\label{o:easybound2} For all  $k \in \bb{N}$, $\bl{\ell}\in \bb{N}^c$, $0< x < 1$ and $\bl{y}\in \bb{R}_{>0}^c$ such that $x+ \sum_{i=1}^c y_i \leq 1$  we have
	$$R(k,\bl{\ell}) \leq x^{-k+1}\prod_{i=1}^{c}y_i^{-\ell_i+1}.$$	
\end{obs}
\begin{proof}
	We will prove this by induction on $\ell$.
	For the induction step, we may assume that $\ell_i \geq 2$ for every $i$, as otherwise the statement trivially holds. (In particular, the base case is trivial.) Let $\bl{\ell^{-,i}}=(\ell_1,\ldots,\ell_{i-1},\ell_{i}-1,\ell_{i+1},\ldots)$ be obtained from $\bl{\ell}$ by subtracting one from $i$-th component. Then we have $$ R(k,\bl{\ell})\leq R(k-1,\bl{\ell}) + \sum_{i=1}^{c}R(k,\bl{\ell^{-,i}})$$
and	applying the induction hypothesis to the terms on the right side we obtain the desired bound. 
\end{proof}

The definitions related to candidates extend as follows.
Let $X$ and $Y$ be two non-empty disjoint subsets of vertices  of a complete graph with edges colored in $c+1$ colors.  We say that $(X,Y)$ is a \emph{$c$-candidate}. For $k \in \bb{N}$ and $\bl{\ell},\bl{t} \in \bb{N}^c$ 
we say that a candidate $(X,Y)$ is \emph{$(k,\bl{\ell},\bl{t})$-good} if $X \cup Y$ contains a red $K_k$ or $X$ contains a monochromatic $K_{t_i}$ colored $B_i$ or $Y$ contains a monochromatic $K_{\ell_i}$ colored $B_i$ for some $i$.  

The following lemma extends \cref{l:easy} to multiple colors. The proof is essentially the same, but we include it for completeness. Recall that  $f_{p}(X,Y)=e_R(X,Y) - p|X||Y|$.

\begin{lem}\label{l:easy2}
	Let  $0<x<p<1$, let $\bl{\theta} \in  \bb{R}_{>0}^c$  be such that $\sum_{i=1}^c \theta_i=1$.  Let $k \in \bb{N}$ and  $\bl{\ell},\bl{t} \in \bb{N}^c$, and let $(X,Y)$ be a $c$-candidate such that 
	\begin{equation}\label{e:easymult}
	f_{p}(X,Y) \geq  (k+t)x^{-k+1}(1-x)^{-\ell+c}(p-x)^{-t+c}\cdot \prod_{i=1}^c\theta_i^{-\ell_i-t_i} 
	\end{equation}	
	then $(X,Y)$ is $(k,\bl{\ell},\bl{t})$-good.
\end{lem}

\begin{proof} 
	We prove the lemma by induction on $k+t$. The cases  when $k=1$ or $t_i=1$ for some $i$, including the base case, are trivial. 
	
For the induction step, note that by  
 \cref{l:FpAvg}, which still applies in the multicolor setting, there exists $v \in X$ such that  $f_p(X,N_R(v) \cap Y) \geq p \cdot f_p(X,Y)$. Let $Y' = N_R(v) \cap Y, X_ {B_i}=  N_{B_i}(v) \cap X$ and $X_R=  N_R(v) \cap X$.
	
	Just as in the proof of \cref{l:easy}, if either $f_p(X_R,Y') \geq \frac{k+t-1}{k+t} x f_p(X,Y)$ or  $f_p(X_{B_i},Y') \geq 
	\theta_i\frac{k+t-1}{k+t}(p-x) f_p(X,Y)$ for some $i$, then applying the induction hypothesis to  $(X_R,Y')$ or  $(X_{B_i},Y')$, respectively, yields the lemma. Thus we assume that none of these inequalities hold and so
	\begin{align*}
		p \cdot f_p(X,Y) &\leq f_p(X,Y') = f_p(X_R,Y') + \sum_{i=1}^c f_p(X_{B_i},Y') + f_p(\{v\},Y') \\ &< \frac{k+t-1}{k+t}pf_p(X,Y) + |Y|,
	\end{align*}
implying
 $\frac{p}{{k+t}}f_p(X,Y) \leq |Y|$.
	Thus, using \cref{o:easybound2},
\begin{align*}
	|Y| &\geq  p x^{-k+1}(1-x)^{-\ell+c}(p-x)^{-t+c}\cdot \prod_{i=1}^c\theta_i^{-\ell_i-t_i}  \geq x^{-k+1}\prod_{i=1}^c\s{ {\underbrace{\theta_i(1-x)}_{=y_i}}}^{-\ell_i+1} \geq  R(k,\bl{\ell}),
\end{align*}	
	implying that $(X,Y)$ is $(k,\bl{\ell},\bl{t})$-good, as desired. 
\end{proof}	

We are ready to prove the multicolor analogue of \cref{t:easy}, which as promised differs from it only by the introduction of the $\Theta(\bl{\ell})$ factor.

\begin{thm}\label{t:easy2}
	For every $\frac{\sqrt{5}-1}{\sqrt{5}+1}< p <1$ and all $k \in \bb{N}$ and  $\bl{\ell} \in \bb{N}^c$ we have $$R(k,\bl{\ell}) \leq 4(k+\ell)\s{\frac{1+\sqrt{5}}{2}p+\frac{1-\sqrt{5}}{2}}^{-k/2}(1-p)^{-\ell} \cdot \Theta(\bl{\ell}).$$
\end{thm}
\begin{proof} Let $\theta_i~=~\ell_i/\ell$. Then $\Theta(\bl{\ell}) = \prod_{i=1}^c \theta_i^{-\ell_i}$. As in \cref{t:easy} the proof is by induction on $\ell$ and we set 
$x = \frac{1+\sqrt{5}}{2}p+\frac{1-\sqrt{5}}{2}$.

Let $Q=4(k+\ell)x^{-k/2}(1-p)^{-\ell}  \prod_{i=1}^c \theta_i^{-\ell_i}$.
Consider a  coloring of the edges of a complete graph on $n \geq \lfloor Q \rfloor$ vertices.  As in the proof of \cref{t:easy}  we may assume that $$|N_{B_i}(v)|\leq  \left \lfloor \theta_i\frac{k+\ell-1}{k+\ell}(1-p)Q \right \rfloor - 1$$ for every vertex $v$ and every $i$, as otherwise the theorem follows by applying the induction hypothesis to $N_{B_i}(v)$. Thus $$|N_R(v)| \geq n-1 - \s{\left \lfloor \frac{k+\ell-1}{k+\ell}(1-p)Q \right\rfloor- 1}$$ for every $v$. As in \cref{t:easy}  we consider a uniformly random partition $(X,Y)$ of the vertices of our graph and following the same steps derive
	$$\bb{E}[f_p(X,Y)] \geq (k+\ell)x^{-k}(1-x)^{-\ell+1}(p-x)^{-\ell+1}\cdot \prod_{i=1}^c\theta_i^{-2\ell_i} .
	$$ Thus $(X,Y)$ is a $(k,\bl{\ell},\bl{\ell})$-good $c$-candidate by \cref{l:easy2}, implying that  our  coloring contains a red $K_k$ or $K_{\ell_i}$ in color $B_i$ for some $i$, as desired.
\end{proof}	

Substituting the same value of $p$ as in \cref{s:easy} into \cref{t:easy2}   yields the following:
\begin{cor}\label{c:easy2}
	For all $k \in \bb{N}$ and  $\bl{\ell} \in \bb{N}^c$ we have $$R(k,\bl{\ell}) \leq 4(k+\ell)\s{\frac{k+2 \ell}{k}}^{k/2}\fs{\left(\sqrt{5}+1\right) (k+2 \ell)}{4\ell}^{\ell}\cdot \Theta(\bl{\ell}).$$
\end{cor}

Just as   in \cref{s:easy}, comparing the bound given in \cref{c:easy2} with the Erd\H{o}s-Szekeres bound we obtain a similar  conclusion
\begin{align*} \frac{R(k,\bl{\ell})}{ES(k,\bl{\ell})} &\leq e^{o(k+\ell)} \fs{(\sqrt{5}+1)(k+2\ell)}{4(k+\ell)}^{\ell}\fs{(k+2\ell)k}{(k+\ell)^2}^{k/2}, 
\end{align*}
and thus still an exponential improvement for $\ell < 0.69 k$. 

\vskip 10pt

Initially, we expected that the results in \cref{s:book} would similarly translate to the multicolor setting, but so far we were unable to overcome the following technical issue. Naively applying a multicolor analogue of \cref{l:BBook} in an extension of \cref{t:bookmain} requires $|X|$ to be sufficiently large compared to $R(k, m_1, t_2,\ldots,t_c) + R(k,t_1,m_2,\ldots,t_c) + \ldots$, where $m_i$ is logarithmic in $k+t$.
Meanwhile, the analogue of \eqref{e:x} gives us a lower bound of the form $|X| \geq (1+\eps)^{k+\ell+t}\Theta(\bl{\ell})$, which is not quite sufficient. 

\begin{rem} After the completion of the initial version of this paper, Balister et al.~\cite{BBCGHMRT} obtained an exponential improvement for the diagonal multicolor Ramsey numbers. 
\end{rem}

\bibliographystyle{halpha}
\bibliography{snorin2}

\appendix

\section{Proof of \cref{lem:numerics}}

We give the explicit formulas and the root-isolation data
used in the numerical check.
For $f=F_{0.08}$, the functions in Lemma~\ref{lem:frontier} are
\begin{align*}
	A(t)&=\frac{t}{1+t}\exp\!\left(-e^{-t}
	\left(-\frac14+\frac{41}{100}t+\frac4{25}t^2
	-\frac2{25}t^3\right)\right),\\
	B(t)&=\frac1{1+t}\exp\!\left(e^{-t}
	\left(\frac{33}{100}t^2+\frac2{25}t^3
	-\frac2{25}t^4\right)\right).
\end{align*}
Thus
\[
a=A(1)=\frac12e^{-0.24/e}=0.4577471843\ldots,
\qquad
b=B(1)=\frac12e^{0.33/e}=0.5645383480\ldots.
\]
Let $\lambda_b$ and $\lambda_a$ be the unique solutions of
$X_1(\lambda_b)=b$ and $X_1(\lambda_a)=a$, respectively.
Directed-rounding interval arithmetic gives the following certificate;
all lower bounds in the last column are rounded down.
\[
\begin{array}{ccl}
	\toprule
	\text{quantity}&\text{interval or location of its minimum}
	&\text{certified lower bound}\\
	\midrule
	\Phi_0(\lambda)/\lambda&(0,1]\ \text{(minimum at $1$)}
	&2.4\cdot10^{-3}\\
	\Phi_1(\lambda)/\lambda,\ X_1\geq b
	&(0,\lambda_b]&1.57\cdot10^{-2}\\
	\Phi_1(\lambda)/\lambda,\ a\leq X_1\leq b
	&[\lambda_b,\lambda_a]&3.5\cdot10^{-4}\\
	\Phi_1(\lambda)/\lambda,\ X_1\leq a
	&[\lambda_a,1]&3.2\cdot10^{-4}\\
	\bottomrule
\end{array}
\]
Here $\Phi_0$ and $\Phi_1$ denote the left sides of
\eqref{eq:base-slack} and \eqref{eq:final-slack}, respectively.
The same interval calculation proves
\[
X_1'(\lambda)<-0.04
\qquad(0<\lambda\leq1),
\]
so the two crossings are unique.  It contains them in
\[
(0.25685117,0.25685119),\qquad
(0.35945039,0.35945042).
\]
For completeness, interval Newton iteration gives the following full
list of interior critical-point enclosures.  It also proves that the
derivative has a fixed nonzero sign on every complementary interval:
\begin{align*}
	(\Phi_0/\lambda)'=0:
	&\quad(0.4067485,0.4067487),\ (0.7929808,0.7929812),\\
	(\Phi_1/\lambda)'=0\quad\text{on }[\lambda_b,\lambda_a]:
	&\quad(0.33831108,0.33831111),\\
	(\Phi_1/\lambda)'=0\quad\text{on }[\lambda_a,1]:
	&\quad(0.4275622,0.4275625),\ (0.6406049,0.6406052),\\[-1mm]
	&\quad(0.8959773,0.8959778).
\end{align*}
There is no critical point on $(0,\lambda_b]$.  Evaluation on these
outward-rounded intervals and at the branch endpoints gives the bounds
in the table.  In particular, the middle-branch minimum is
$0.0003525554\ldots$, and the last-branch minimum is
$0.0003238078\ldots$.  No numerical differentiation is needed: on
the branches $X_1=B(t)$, $X_1=A(t)$, and $a\leq X_1\leq b$, respectively,
implicit differentiation gives
\[
(\log Y_1)'=-\frac1t(\log X_1)',\qquad
(\log Y_1)'=-t(\log X_1)',\qquad
(\log Y_1)'=-(\log X_1)'.
\]
At the removable endpoint $\lambda=0$ one has
\begin{align*}
	\lim_{\lambda\downarrow0}\frac{\Phi_0(\lambda)}\lambda
	&=\frac34-\frac{1+e^{1/4}}2
	+\frac12\log(1+e^{1/4})=0.0209569\ldots,\\
	\lim_{\lambda\downarrow0}\frac{\Phi_1(\lambda)}\lambda
	&=\frac78-\frac{1+e^{1/4}}2
	+\frac12\log(1+e^{1/4})=0.1459569\ldots.
\end{align*}
This completes the verification of \eqref{eq:base-slack} and
\eqref{eq:final-slack}.

\end{document}